\documentclass[a4paper,11pt,oneside]{amsart}

\usepackage{geometry}
\usepackage{graphicx, amsfonts, amssymb}
\usepackage{mathrsfs, amsmath, amsthm}
\usepackage{verbatim}
\usepackage{relsize}
\usepackage{enumerate}
\usepackage{hyperref}
\usepackage{version}
\usepackage{blindtext}
\usepackage{mathtools}
\usepackage{color}

\newtheorem{theorem}{Theorem}[section]

\newtheorem{definition}[theorem]{Definition}
\newtheorem{Lemm}[theorem]{Lemma}
\newtheorem{prop}[theorem]{Proposition}
\newtheorem{corollary}[theorem]{Corollary}

\theoremstyle{remark}

\newcommand{\m}{\mathbb}
\newcommand{\mm}{\mathcal}

\def\CB{\color{black} }

\DeclareMathOperator{\Hol}{Hol}

\DeclareMathOperator{\BMOA}{BMOA}
\DeclareMathOperator{\Mer}{Mer}
\DeclareMathOperator{\loc}{loc}
\hypersetup{
	colorlinks=true,
	linkcolor=blue,
	filecolor=magenta,      
	urlcolor=cyan,
	pdftitle={Meromorphic Optimal Domain of Integral Operators},
	pdfpagemode=FullScreen,
}

\title{Meromorphic Optimal Domain of Integral Operators}

\author{C. Bellavita}
\email{carlobellavita@ub.edu}
\address{Departamento de Matem\'atica i Inform\'atica, Universitat de Barcelona,
Gran Via 585, 08007 Barcelona, Spain}

\author{A. Belli}
\email{anilbelli@math.auth.gr}
\address{Department of Mathematics, Aristotle University of Thessaloniki, 54124, Greece}

\author{G. Nikolaidis}
\email{nikolaidg@math.auth.gr}
\address{Department of Mathematics, Aristotle University of Thessaloniki, 54124, Greece}

\author{G. Stylogiannis}
\email{stylog@math.auth.gr}
\address{Department of Mathematics, Aristotle University of Thessaloniki, 54124, Greece}

\date{\today}

\begin{document}

\begin{abstract}
    For $g\in\BMOA$, we introduce the meromorphic optimal domain $(T_g,H^p)$, i.e. the space containing the meromorphic functions that are mapped under the action of the generalized Volterra operator $T_g$ into the Hardy space $H^p$. We investigate its properties and characterize for which $g_1,g_2 \in \BMOA$ the corresponding meromorphic optimal domains coincide. 
    This investigation contributes to a more comprehensive understanding of the holomorphic optimal domain of $T_g$ in $H^p$.
\end{abstract}
\maketitle

\section{Introduction}
Let $\m{D}$ be the unit disc of the complex plane, $\m{T}$ be its boundary, and  $\Hol(\m{D})$ the space of analytic functions defined in $\m{D}$. The (shifted) Ces\'aro summation operator is defined as
$$
C(f)(z)=\int_{0}^{z}\frac{f(\zeta)}{1-\zeta}\,d\zeta\, .
$$
For a fixed $g\in \Hol(\m{D})$,  we consider the generalized Volterra operator, acting on the space of holomorphic functions, as 
\begin{equation*}\label{equation for Tg}
    T_g(f)(z)=\int_{0}^{z}f(\zeta)g'(\zeta)\,d\zeta\,,\qquad z\in\m{D},\quad f\in\Hol(\m{D})\,.
\end{equation*}
The motivation for studying this operator arises from the fact that, as $g$ varies within $\Hol(\m{D})$, $T_g$ represents significant classical operators. 
For instance, if $g(z)=z$, then $T_g$ is the classical Volterra operator, while, choosing $g(z)=-\log(1-z)$, then $T_g$ coincides with the Ces\'aro summation operator. 

The primary focus in studying $T_g$ operators is their actions on spaces of analytic functions, with particular emphasis on the Hardy spaces of the unit disc. Let $0<p<\infty$, the Hardy space $H^p$ consists of all functions $f\in \Hol(\mathbb{D})$ such that

$$   \|f\|^p_{H^p}=\sup_{0\leq r<1}\int_{0}^{2\pi}|f(re^{it})|^p\,\frac{dt}{2\pi}<\infty\,.$$
	For $p=\infty$, $H^{\infty}$ denotes the space of bounded analytic functions in the unit disc, equipped with the supremum norm.

\par Ch. Pommerenke \cite{pommerenke1977schlichte} proved that the space of symbols $g$ for which $T_g$ is bounded on the Hilbert space $H^2$ is $\BMOA$. In continuation, A. Aleman and A. Siskakis \cite{aleman1995integral} showed that the space of symbols $g$ for which $T_g$ is bounded on the Hardy space $H^p$, for $1\leq p<\infty$, is still $\BMOA$. The space of analytic functions of bounded mean oscillation, $\BMOA$, consists of all $f \in H^2$ such that 
$$ 
\|f\|_*:=\sup_{a\in \mathbb{D}} \| f \circ \varphi_a - f(a) \|_{H^2} < \infty,
\text{ where } \varphi_a(z) = \frac{z-a}{1-\overline{a}z}\, .
$$
 Subsequently, the study of the boundedness and the properties \CB of $T_g$ operators on various spaces of analytic functions attracted a lot of attention, see, for example, \cite{aleman1997integration}, \cite{wu2006areabergman}, \cite{pelaez2014weightedmemoir}, \cite{TanausuGalanopoulostentspaces}, \cite{Chalmoukis2024} and \cite{bellavita2024}. 
 
 For fixed $g\in\BMOA$, even though $T_g(H^p)\subset H^p$, there exists function $f$ not belonging to $H^p$, for which $T_g(f)$ is still in $H^p$, see \cite[Theorem 1]{Bellavitaoptimaldom}. In \cite{Curbera2011} and \cite{curbera2012}, G. Curbera and W. Ricker introduced the optimal domain for the Ces\'aro operator $[C,H^p]$, consisting of all the analytic functions which are mapped by $C$ into $H^p$. In \cite{Bellavitaoptimaldom}, the authors defined and studied the optimal domain for integral operator $T_g$ 
 $$
 [T_g,H^p]=\{f\in \Hol(\m{D})\colon T_g(f) \in H^p\}\, ,
 $$
 It consists of all the analytic functions mapped by $T_g$ into $H^p$.
 
 \par It appears that the main tool needed to study $[C,H^p]$, is its equivalent description \cite[Proposition 3.2]{Curbera2011}
 \begin{equation}\label{E:exp[C,H]}
 [C,H^p]=\{f\in \Hol(\m{D})\colon f(z)=(1-z)h'(z)\text{ for some } h\in H^p\}\, ,
 \end{equation} 
 from which it follows that the Ces\'aro operator is surjective from $[C,H^p]$ onto $H^p_{0}$.
 In \cite{Bellavitaoptimaldom}, the authors observed that such property does not hold for $[T_g,H^p]$ for a general symbol $g\in\BMOA$, since $g'$ may have zeros inside the unit disc. Hence, it is clear that we need to relax the condition $f\in \Hol(\m{D})$. We denote as $\Mer(\m{D})$ the set of meromorphic functions defined in the unit disc $\m{D}$.
 \begin{definition}\label{Meromorphic Optimal Domain}
     Let $g\in\BMOA$ and $1\leq p<\infty$. The meromorphic optimal domain $(T_g,H^p)$, consists of all the functions $f\in \Mer(\m{D})$, such that $f g'\in \Hol(\m{D})$ and $T_g(f)\in H^p.$
 \end{definition}
 The condition $fg'\in \Hol(\m{D})$ is necessary to make $T_g(f)$ well defined. Also, the same condition implies that the poles of $f$ belong to the zero set of $g'$ with at most the same multiplicity.

\vspace{11 pt}
\par
 Characterizing all the functions $h \in \BMOA$ such that $(T_h,H^p)= (T_g,H^p)$ is an interesting \CB problem. To address it, we define and study the following space.
  \begin{definition}\label{Wg definition}
     Let $1\leq p<\infty$ and $g\in\BMOA$. We consider the spaces
     \begin{align*}
     W_g^p&:=\left\lbrace h \in \BMOA \colon (T_h,H^p)\supseteq (T_g,H^p)\right\rbrace\, ,\\
     V_g^p&:=\left\lbrace h \in \BMOA \colon [T_h,H^p]\supseteq [T_g,H^p]\right\rbrace\, .
     \end{align*}
 \end{definition}
 
\par The primary aim of this article is to examine the properties of $(T_g,H^p)$, $W^p_g$ and $V^p_g$ .Through this study, we enhance our knowledge of the holomorphic optimal domain $[T_g,H^p]$ and we uncover a noteworthy connection between $W^p_g$, $V^p_g$ and $T_g(H^\infty)\subset \BMOA$. Finally, we discuss given $1\leq p<\infty$ and $g_1,g_2\in\BMOA$, whether
$(T_{g_1},H^{p})=(T_{g_2},H^p)$ and $[T_{g_1},H^{p}]=[T_{g_2},H^p]$ hold. 

 \par
 As it is customary, for real-valued functions $f,\, g$, we write $f\lesssim g$, if there exists a positive constant $C$ (which may be different in each occurrence) independent of the arguments of $f,\,g$ such that $ f\leq Cg$. The notation $f\gtrsim g$ can be understood analogously. If both $f\lesssim g$ and $f\gtrsim g$ hold simultaneously, then we write $f\sim g$. We use the symbol $\mathcal{Z}(h)$ to denote the zero set of $h\in\Hol(\m{D})$, where the elements of $\mathcal{Z}(h)$ are taken according to their multiplicities.

\vspace{22 pt}
 \section{Main Results}
 For $g\in\BMOA$, we note that $H^p\subset (T_g,H^p)$. By Lemma \ref{Lemma Meromorphic constant}, we have that $(T_g,H^p)$ is the space of meromorphic functions only when $g$ is constant. So, for non constant symbol $g$, we endow $(T_g,H^p)$ with the norm
 \begin{equation}\label{E:norm T_G}
   \|f\|_{(T_g,H^p)}:= \|T_g(f)\|_{H^p}\, .
 \end{equation}
In order to study the properties of the meromorphic optimal domains, we remark that $(T_g,H^p)$ has the following description
\begin{equation}\label{E:relation meromorphic}
    (T_g,H^p)=\biggl\{ f\in \Mer(\m{D})\colon f(z)=\frac{h'(z)}{g'(z)}  \text{ for some } h\in H^p\biggr\}
\end{equation}
for $g \in \BMOA$, non-constant. Since the proof of this result follows the line of \cite[Proposition 3.2]{Curbera2011}, we omit its proof. The main results of our study about properties of $(T_g,H^p)$ are summarized in the following theorem.
\begin{theorem}\label{Theorem properties of meromorphic optimal domain}
    Let $g\in\BMOA$, non-constant and $1\leq p<\infty$. The following properties hold:
    \begin{itemize}
        \item[(a)] The pair $\left( (T_g,H^p), \|\cdot\|_{(T_g,H^p)}\right)$ is a Banach space of meromorphic functions with bounded point evaluation functional at every $z_0 \in \mathbb D\setminus \mathcal{Z}(g')$.
        \item[(b)] If $1\leq p_1<p_2<\infty$, we have that
        $$(T_g,H^{p_2})\subsetneqq(T_g,H^{p_1})\,.$$
        \item[(c)] Let $\mm{M}(T_g,H^p)$ be the space of multipliers of $(T_g,H^p)$. Then, $\mm{M}(T_g,H^p)=H^{\infty}$.
        \item[(d)] The space $(T_g,H^p)$ is always separable.
        \item[(e)] Let $1<p<\infty$. The dual of $(T_g,H^p)$ can be identified with $(T_g,H^{q})$ with $\frac{1}{p}+\frac{1}{q}=1$, under the pairing
        $$
        \left\langle f, k\right\rangle=\left\langle T_g(f), T_g(k)\right\rangle_{L^2(\m{T})}\qquad f\in (T_g,H^p),\,k\in (T_g,H^q)\,.
        $$
        \item[(f)] Let $1<p<\infty$.  $(T_g,H^p)$ is the $\left(1-1/p,p\right)$-interpolation space between $(T_g,H^1)$ and $(T_g,H^{\infty})$, that is, by using the notation of \cite[Ch. 5.1, Def. 1.7]{bennett1988interpolation},
        $$
        ((T_g,H^1),(T_g,H^{\infty}))_{1-\frac{1}{p},p}=(T_g,H^p)\, .
        $$
    \end{itemize}
\end{theorem}
Theorem \ref{Theorem properties of meromorphic optimal domain} $(d)$ implies that every holomorphic optimal domain $[T_g,H^p]$ is separable as well, see Corollary \ref{Corollary of standard optimal domain separability}, even if the polynomials do not always form dense subsets of $[T_g,H^2]$, see \cite{Bellavitaoptimaldom}. \CB

   \par We turn  now \CB our attention to $W_g^p$. We are interested in describing for which $g_1,g_2 \in \BMOA$ the corresponding meromorphic optimal domains are equal, that is $(T_{g_1},H^p)=(T_{g_2},H^p)$.
   From the definition of $W_{g_i}^p$, it is clear that $(T_{g_1},H^p)=(T_{g_2},H^p)$ if and only in $g_1 \in W^p_{g_2}$ and $g_2 \in W^p_{g_1}$. For this reason, having an explicit description of $W^p_g$ is fundamental.
    \begin{theorem}\label{Thereom characterization of Wg with Tg(Hoo)}
      Let $g\in\BMOA$ and $1\leq p<\infty$. $h\in W_g^p$ if and only if there exists a function $f\in T_g(H^{\infty})$ such that $h=f+h(0)$. In other words,
      $$W_g^p=T_g(H^{\infty})+\m{C}\,.$$
  \end{theorem}
  Lemma \ref{Lemma Meromorphic constant} implies that all the constant functions belong to $W_g^p$ for every $g \in \BMOA$. We notice also that $W_g^p$ is independent of $p$ and, hence, we can omit the superscript $p$.  Based on Theorem \ref{Thereom characterization of Wg with Tg(Hoo)}, we answer when two meromorphic optimal domains coincide.
\begin{theorem}\label{T:equivalence meromorphic opttimal}
 For $1\leq p<\infty$, let $g_1,g_2\in\BMOA$, non-constant. The following conditions are equivalent:
 \begin{itemize}
     \item[a)] $(T_{g_1},H^p)=(T_{g_2},H^p)$.
     \item[b)] There exist two functions $k_1,k_2\in H^{\infty}$ such that 
     $$
g_1=T_{g_2}(k_1)+g_1(0) \text{ and } g_2=T_{g_1}(k_2)+g_2(0)\,.
$$
 \end{itemize}
 Moreover, 
 we have that $k_2=1/k_1$\,.
\end{theorem}  
  
Afterward, we investigate the structural properties of $W_g$. To this end, we endow $W_g$ with the following norm: for $h \in W_g$ and $h-h(0)=T_g(k)$ with $k \in H^\infty$,
  $$
\|h\|_{W_g}=\|T_g(k)\|_{*}+|h(0)|+\|k\|_{H^\infty}\, .
  $$ 
The following theorem summarizes our study of the properties of $W_g$.  \begin{theorem}\label{Theorem properties of Wg}
      Let $g\in\BMOA$, non-constant and $1< p<\infty$. \CB The following properties hold:
      \begin{itemize}
          \item[(a)]The pair $\left( W_g, \|\cdot\|_{W_g}\right)$ is a linear subset of $\BMOA$ and it is a Banach space.
          \item[(b)] We have that $$
    \bigcap_{h \in W_g\colon (T_h,H^p)\supsetneq (T_g,H^p)} (T_h,H^p)=(T_g,H^p)\,  .
$$
\item[(c)] Let $\mm{M}(W_g)$ denote the set of multipliers of $W_g$. Then $\phi\in\mm{M}(W_g)$ if and only if $\phi\in H^{\infty}$ and the integral operator $T_\phi$ is a bounded operator on $W_g.$
      \end{itemize}
  \end{theorem}

  For some $g \in \BMOA$, the space $W_g$ equipped with the norm $\|\cdot\|_{*}$ is not a closed subspace of $\BMOA$. For this reason, in $(a)$ it is fundamental to consider $\|\cdot\|_{W_g}$ as its rightful norm. 
 
  \par The study of the meromorphic optimal domains yields results also for the holomorphic optimal domains. For instance, let $U_{\loc}(\m{D})$ be the set of locally univalent functions in $\m{D}.$ It is well-known that $g\in U_{\loc}(\m{D})$ if and only if $g'$ never vanishes in $\m{D}$.  In this case, it is easy to check that
  $$
  (T_g,H^p)=[T_g,H^p]\, .
  $$
   Consequently, when $g \in U_{\loc}(\m{D})$, we obtain an explicit description of $V_g^p$, see Proposition \ref{Proposition characterisation Vg g loc univ}, which reformulates the idea used in the proof of Theorem \ref{Thereom characterization of Wg with Tg(Hoo)}. We also provide some inclusion \CB concerning $V_g$ when $g'$ is factorized into an $H^\infty$ function times the derivative of a $\BMOA$ function, \CB Proposition \ref{proposition V_g when g in nevanlinna}.
  \par The final property we study in this article is whether two holomorphic optimal domains coincide. However, this question is extremely challenging. To even out the situation, we impose certain natural conditions on the symbols $g_i$.
 \begin{itemize}
     \item[(a)] We deal with the case when $g_1,g_2 \in \BMOA \cap U_{\loc}(\m{D})$. We obtain the answer (Proposition \ref{prop univalent}) by reformulating the idea used in the proof of Theorem \ref{Thereom characterization of Wg with Tg(Hoo)}.
     
     \item[(b)] In Proposition \ref{Theorem optimal Domain polynomial}, we consider the case when $g$ is an analytic polynomial and we describe the connection between  $[T_z,H^p]$ and $[T_g,H^p]$.
     \item[(c)] In Proposition \ref{p=2 optimaldomain question}, we consider the Hilbert case $p=2$ and we reformulate the problem in terms of multipliers bounded from below acting on non-radial weighted Bergman spaces. 
 \end{itemize}

  \par The rest of the article is divided into three sections: in Section 3, we focus our attention on the meromorphic optimal domain and prove Theorem  \ref{Theorem properties of meromorphic optimal domain}. In Section 4, we describe the space $W_g$ and provide the proof of Theorem \ref{Thereom characterization of Wg with Tg(Hoo)}, \ref{T:equivalence meromorphic opttimal} and \ref{Theorem properties of Wg}. Finally, in Section 5, we briefly describe the structure of $V^p_g$ and show when two holomorphic optimal domains coincide in some particular situations.
  
  \CB
  
 \vspace{22 pt}
\section{Linear space structure of \texorpdfstring{$(T_g,H^p)$}{}.}
First of all, we describe $(T_g,H^p)$ when $g$ is a constant function. \CB 
\begin{Lemm}\label{Lemma Meromorphic constant}
 Let $g \in \BMOA$ and $1\leq p<\infty$. Then $(T_g,H^p)=\Mer(\m{D})$ if and only if $g$ is constant.
 \end{Lemm}
\begin{proof}
Let $f$ be a meromorphic function. If $g$ is constant, then $fg'$ is the zero function and, consequently, it is analytic. Since in this case $T_g(f)\equiv 0$, the meromorphic optimal domain coincides with $\Mer(\m{D})$. 

We prove now the necessity and we assume that $(T_g,H^p)=\Mer(\m{D})$. Let $z_0\in\m{D}$. The function 
$$f(z)=\frac{1}{z-z_0}\qquad z\in\m{D}$$
is the meromorphic function, with only a pole of order one at $z_0$. Our assumption implies that $f\in (T_g,H^p)$. Consequently,
$$\frac{g'(z)}{z-z_0}\qquad z\in\m{D}$$
is analytic in the unit disc. This fact implies that $g'$ has a root of at least multiplicity one at $z_0$. Since $z_0$ is arbitrarily chosen, it implies that $g'$ is zero everywhere on $\m{D}$ and consequently $g$ needs to be constant.\\
\end{proof}

Due to the linearity of the operator $T_g$, $(T_g,H^p)$ is a vector space of meromorphic functions in the unit disc. The quantity $\|\cdot \|_{(T_g,H^p)}$, defined in \eqref{E:norm T_G}, introduces a seminorm, as it is subadditive and scalar multiplicative.
\begin{Lemm}\label{Tg injectivity true norm Lemma}
    Let $g\in\BMOA$ and non-constant. Then $\|f\|_{(T_g,H^p)}=0$ if and only if $f=0$. Hence, $\|\cdot\|_{(T_g,H^p)}$ is a norm in $(T_g,H^p)$. 
\end{Lemm}
\begin{proof}
It is clear that $T_g(0)=0$. Let $f\in \Mer(\m{D})$ such that $fg'\in \Hol(\m{D})$ and $\|T_g(f)\|_{H^p}=0$. Then, $T_g(f)(z)=0$,  for every $z\in \m{D}$. Differentiating, we obtain that
$$
g'(z)f(z)=0\,\, \text{ for all } z\in\m{D}\, .    
$$
Since $g' \neq 0$, then $f(z)=0$ for every $z \in \m{D}\setminus \mathcal{Z}(g').$ As $\mathcal{Z}(g')$ is at most countable, it means that the $\displaystyle{\lim_{z\to w}f(z)=0}$ for every $w\in \mathcal{Z}(g')$. This implies that $f$ has no poles and that $f$ is zero almost everywhere on $\m{D}$. Consequently, $f=0.$\\
\end{proof}
Before showing the proof of Theorem \ref{Theorem properties of meromorphic optimal domain}, we prove that the point evaluation functionals are bounded. 
\begin{Lemm}\label{Lemma point evaluation functional}
 Let $g\in\BMOA$, be non-constant. Let $z_0\in\m{D}$ and $\lambda_{z_0}$ denote the point evaluation functional at $z_0.$ The functional $\lambda_{z_0}$ is bounded in $(T_g,H^p)$ if and only if $z_0 \in \mathbb D\setminus \mathcal{Z}(g')$.
\end{Lemm}
\begin{proof}
Let $z_0\in \mathbb D\setminus \mathcal{Z}(g')$ and $f\in (T_g,H^p)$. By standard argument involving the point evaluation functionals for the derivatives of $H^p$ functions, we have that
    $$
    |T_g(f)'(z_0)|{\lesssim \frac{1}{(1-|z_0|)^{\frac{1}{p}+1}}\|T_g(f)\|_{H^p}\, .}
    $$
    Since $\|f\|_{(T_g,H^p)}= \|T_g(f)\|_{H^p}$, we note that
    $$|f(z_0)\cdot g'(z_0)|{\lesssim \frac{1}{(1-|z_0|)^{\frac{1}{p}+1}}\|f\|_{(T_g,H^p)}}\, .
    $$
    and
    $$
    |f(z_0)|{\lesssim \frac{1}{|g'(z_0)|(1-|z_0|)^{\frac{1}{p}+1}}\|f\|_{(T_g,H^p)}\, .}
    $$
On the other hand, let $z_0\in \mathcal{Z}(g')$. As $1/g'\in (T_g,H^p)$, and $\lim_{z\to z_0}{1}/{g'(z)}=\infty$, we conclude that $\lambda_{z_0}$ is unbounded.\\
\end{proof}
\begin{proof}[Proof of Theorem \ref{Theorem properties of meromorphic optimal domain}] (a)
Due to Lemma \ref{Tg injectivity true norm Lemma} and Lemma \ref{Lemma point evaluation functional}, we have only \CB to prove that $\left( (T_g,H^p), \|\cdot\|_{(T_g,H^p)}\right)$ is a Banach space.
    Let $\{f_n\}_n$ be a Cauchy sequence in $(T_g,H^p)$, that is
    \begin{equation}\label{equi2}
        \|f_n-f_m\|_{(T_g,H^p)}=\|T_g(f_n-f_m)\|_{H^p}\xrightarrow[]{n,m\to+\infty}0\, .
    \end{equation}
    Since $H^p$ is a Banach space, \eqref{equi2} implies the existence of a function $k\in H^p$, such that
    $$
    T_g(f_n)\xrightarrow[n\to+\infty]{\|\,\|_{H^p}} k\, .
    $$
    Since the point evaluation for the derivative of the Hardy spaces are uniformly bounded on compact sets, we have that 
    $$
    f_n(z)g'(z)=T_g(f_n)'(z)\xrightarrow[n\to+\infty]{}k'(z)\qquad \text{ for every }z \in \mathbb D\, .
    $$
 Moreover, due to \eqref{E:relation meromorphic}, the meromorphic function $f=k'/g'$ belongs to $(T_g,H^p)$ and it is the limit point of $\{f_n\}_n$ since \CB
$$
\|f_n- f\|_{(T_g,H^p)}=\|T_g(f_n)-k\|_{H^p}\xrightarrow[n\to\infty]{} 0\, .
$$

     (b) The inclusion $(T_g,H^{p_2})\subseteq (T_g,H^{p_1})$ follows from the fact that $H^{p_2}\subset H^{p_1}$.
     To prove that  this inclusion is strict, we consider the function
$$
\phi(z)=\frac{1}{(1-z)^{\frac{1}{p_2}+1}}\frac{1}{g'(z)}\, .
$$
As every zero of $g'$ in $\mathbb{D}$ is of finite order, $\phi\in \Mer(\m{D}).$ Moreover, $\phi g'\in \Hol(\m{D})$ and
$$
T_g(\phi)(z)=-p_2\left(\frac{1}{(1-z)^{1/p_2}}-1\right)\, .
$$
Standard arguments verify that $T_g(\phi)\in H^{p_1}\setminus H^{p_2}$ and the statement is proved.
 
 (c) We recall that for $g,F\in\Hol(\m{D})$, 
$$
S_g(F)(z)=\int_{0}^{z}g(\zeta)F'(\zeta)\,d\zeta
$$
is the $T_g$-companion operator.
 
    Let $f\in (T_g,H^p)$ and $k \in H^\infty$. We want to show that $kf\in (T_g,H^p).$ By \eqref{E:relation meromorphic}, let $h\in H^p$ such that $f=h'/g'$. Then 
    $$ 
    k(z)f(z)g'(z)=k(z)h'(z)\qquad \forall z\in\m{D}\, .
    $$
    As $fg'\in \Hol(\m{D})$, $kf \in \Mer(\m{D})$ and $kfg'\in \Hol(\m{D})$. \CB Moreover, we have that
    $$
    T_g(kf)(z)=S_{k}(h)(z)\qquad z\in \m{D}\, .
    $$
   By \cite[Theorem 2.2]{Anderson2011}, \CB $S_k(h)\in H^p$ and consequently $T_g(kf)\in H^p$ as well.

    On the other hand, let $\phi\in \Mer(\mathbb{D})$, such that $\phi f\in (T_g,H^p)$ for every $ f\in (T_g,H^p)$. 
    Since the function $1/{g'} \in (T_g,H^p)$, the definition of $(T_g,H^p)$ implies that $\phi \in \Hol(\m{D})$. 
    By the result of Lemma \ref{Lemma point evaluation functional}, when $z_0 \in \m{D}\setminus \mathcal{Z}(g')$, $\lambda_{z_0}$ is a bounded linear functional and since $(T_g,H^p)$ contains constant functions, $\|\lambda_{z_0}\|\neq 0$. For $f \in (T_g,H^p)$ with $\|f\|_{(T_g,H^p)} = 1$, we have
$$
|\phi(z_0)f(z_0)| = |\lambda_{z_0}(M_\phi(f))| \leq \|\lambda_{z_0}\|\|M_\phi\|\, .
$$
Thus, taking the supremum on $f \in  (T_g,H^p)$ with $\|f\|_{(T_g,H^p)} = 1$, we have that
$$
|\phi(z_0)|\leq \|M_\phi\|\, .
$$
As $z_0$ is arbitrary, the above relation holds in $\m{D}\setminus \mathcal{Z}(g')$, and, since $\phi\in \Hol(\m{D})$ and $\mathcal{Z}(g')$ is countable, by taking limits, we have that
for every $z \in \m{D}$,
$$
|\phi(z)|\leq \|M_{\phi}\|\, ,
$$
from which the statement follows.

(d) Let $f\in (T_g,H^p)$. By the density of the polynomials in $H^p$, see \cite{duren1970theory}, \CB there is a sequence of polynomials $\{P_n\}_n$ with $P_n(0)=0$ for each $n$, such that
    $$\|T_g(f)-P_n\|_{H^p}\xrightarrow[n\to+\infty]{}0\,.$$
   We consider the sequence 
    $$
    f_n(z)=\frac{P_n'(z)}{g'(z)}\quad \text{ with }z\in\m{D}\, .
    $$
   Due to \eqref{E:relation meromorphic}, $f_n\in (T_g,H^p)$ and
    $$\|f-f_n\|_{(T_g,H^p)}=\|T_g(f)-P_n\|_{H^p}\xrightarrow[n\to+\infty]{}0\, .$$
    Therefore the countable set \CB
    $$\biggl\{f\in \Mer(\m{D})\colon f(z)=\frac{P_n'(z)}{g'(z)}\, \, \text{ with } P_n\text{ {polynomial with rational coefficients}}\biggr\}$$
   is dense in $(T_g,H^p)$.
\par (e) We note that it suffices to prove that $(T_g,H^p)$ is isometrically isomorphic to $H^p_0=\{f\in H^p\colon f(0)=0\}.$ We properly describe the isomorphism between these spaces. Let 
$$
\mm{F}\colon (T_g,H^p)\to H^p_0\colon f\to T_g(f)\, .
$$
From the definition of $(T_g,H^p)$, $\mm{F}$ is well defined. As $T_g$ is injective in $\text{Hol}(\mathbb{D})$, we have that $\mm{F}$ is injective. For the surjectivity, let $h\in H^p_0$. Then, due to \eqref{E:relation meromorphic}, the function
$$
f(z)=\frac{h'(z)}{g'(z)}\
$$
belongs to $(T_g,H^p)$ and  $\mm{F}(f)=T_g(f)=h$.
\par Finally, the isometry can be easily checked since
$$
\|\mm{F}(f)\|_{H^p}=\|T_g(f)\|_{H^p}=\|f\|_{(T_g,H^p)}\, .
$$
Given the isometric isomorphism between $\left( H^p_0\right)^*$ and $H^q_0$ with the duality product defined by the $L^2(\m{T})$ inner product, the statement follows.\CB

\par (f) To prove this statement, we introduce for the convenience of the reader some further notions. Let $X_0,X_1$ be a compatible pair of Banach spaces. For given $f\in X_0+X_1$ and $t>0$, the Peetre's $K$-functional is defined as
$$K(f,t,X_0,X_1)=\inf\{\|f_0\|_{X_0}+t\|f_1\|_{X_1}\colon f=f_0+f_1\,,\, f_0\in X_0,\, f_1\in X_1\}\,.$$
For such pair, the interpolation space defined via the $K$-method, is denoted by $(X_0,X_1)_{\theta,q}$ and it contains all the functions $f\in X_0+X_1$ such that
$$\|f\|_{\theta,q}=\left(\int_{0}^{\infty}t^{-\theta}K(f,t,X_0,X_1)\frac{dt}{t}\right)^{\frac{1}{q}}<\infty\qquad 0<\theta<1, 1\leq q<\infty\,  .$$
 \par We turn our attention to the specific case, where $X_0=(T_g,H^{\infty})$, $X_1=(T_g,H^1)\,,$ $\theta=1-\frac{1}{p}$ and $q=p$. To be precise, we first observe that $(T_g,H^\infty)$ is not empty since it always contains $1/g'$. 
 
As in (e), we know that $({T}_g,H^1)$ and $({T}_g,H^{\infty})$ are isometrically isomorphic to $H^1_0$ and $H^\infty_0$ respectively. Consequently, we fix $f\in ({T}_g,H^{1})=({T}_g,H^1)+({T}_g,H^\infty)$ with $f = h_1 + h_2$, where $h_1 \in (T_g, H^1)$ and $h_2 \in (T_g, H^\infty)$. We note that
$$
T_g(f) = H_1 + H_2
$$
with $H_1=T_g(h_1) \in H^1_0$ and $H_2=T_g(h_2) \in H^\infty_0$. 
The isometry between
$(T_g, H^p)$ and $H^p_0$ implies that
\begin{align}
K(f, t; (T_g, H^1), (T_g, H^\infty))&=\inf\{\|h_1\|_{(T_g,H^1)}+t\|h_2\|_{(T_g,H^{\infty})}\colon f=h_1+h_2\}\nonumber\\
&=\inf\{\|H_1\|_{H^1}+t\|H_2\|_{H^\infty}\colon T_g(f)=H_1+H_2\}\nonumber\\
&= K(T_g(f), t; H^1_0, H^\infty_0),\quad t>0 \label{equation K functional}\, .
\end{align}
We recall that  
\begin{equation}\label{interpolation hp0}
    (H^1_0,H^{\infty}_0)_{1-\frac{1}{p},p}=H^p_0\, ,
\end{equation}
since $H^p_0=Q(L^p)$ where $Q$ is a Calderon-Zygmund convolution operator, see \cite{Kislyakovshuinterpolation} and \cite{bourgain1992some}.
Employing \eqref{interpolation hp0} in combination with \eqref{equation K functional}, we obtain that $((T_g,H^1),(T_g,H^{\infty}))_{1-\frac{1}{p},p}$ contains all the meromorphic functions $f$ for which $T_g(f)\in H^p_0$. In other words, 
$$((T_g,H^1),(T_g,H^{\infty}))_{1-\frac{1}{p},p}=(T_g,H^p)\, .$$
\end{proof}

We point out that since the meromorphic optimal domains are separable, then also the holomorphic optimal domains are separable.

\begin{corollary}\label{Corollary of standard optimal domain separability}
    Let $1\leq p<\infty$, $g\in\BMOA$ and non-constant. Then $[T_g,H^p]$ is separable.
\end{corollary}
\begin{proof}
    We recall that from \cite[Theorem 1]{Bellavitaoptimaldom},  $[T_g,H^p]$ is a closed subspace of the Banach space $(T_g,H^p)$. Moreover, Theorem \ref{Theorem properties of meromorphic optimal domain} part d), implies that $(T_g,H^p)$ is separable. Since the separability is preserved in metric subspaces, the conclusion follows.\\
\end{proof}
Since  $[T_g,H^2]$ is isomorphic to the non-radial Bergman space $A^2(|g'(z)|^2(1-|z|^2))$, see \cite{Bellavitaoptimaldom}, these Bergman spaces are separable as well. Up to the authors knowledge, this property of the Bergman spaces is new and non-trivial, since, for some $g\in \BMOA$, the polynomials are not dense in  $A^2(|g'(z)|^2(1-|z|^2))$, see \cite[Proposition 4]{Bellavitaoptimaldom}.

\vspace{22 pt}
\section{The space \texorpdfstring{$W_g$}{}.}
In this section we discuss when $(T_{g_1},H^p)=(T_{g_2},H^p)$. To this end, the space $W_g^p$ plays an important role, as we shall verify below. 
\begin{proof}[Proof of Theorem \ref{Thereom characterization of Wg with Tg(Hoo)}]
    The constant functions belong to $W^p_g$, due to Lemma \ref{Lemma Meromorphic constant}. Let $f-f(0)=T_g(k)$ with $k\in H^{\infty}.$ Then, for every $z \in \m{D}$, 
    $$
    f'(z)=k(z)g'(z)\, .
    $$
   We want to show that $f\in W_g^p$, that is, $(T_f,H^p)\supseteq (T_g,H^p).$ Let $G\in (T_g,H^p)$. Then, due to Theorem \ref{Theorem properties of meromorphic optimal domain}.c, $kG\in (T_g,H^p)$ and 
   $$
   T_f(G)=T_g(kG)\in H^p\, ,
   $$
   which proves the first inclusion.
   \par On the other hand, let $f \in W^p_g$. Since $1/g' \in (T_g,H^p)$ and $(T_g,H^p)\subseteq (T_f,H^p),$ we have $1/g' \in (T_f,H^p)$. Hence, due to the definition of meromorphic optimal domain,  
   $$
   \frac{f'}{g'} \in \Hol(\m{D})\, .
   $$
    We consider now any $G\in (T_g,H^p)\subseteq (T_f,H^p)$. Then $T_f(G)\in H^p$ and
    $$T_g\left(\frac{f'}{g'}G\right)=T_f(G)\in H^p.$$
    Consequently, $f'/g'$ is a multiplier of $(T_g,H^p)$ and, due to Theorem \ref{Theorem properties of meromorphic optimal domain}.c, we have $f'/g'\in H^{\infty}.$ This proves that $f-f(0)\in T_g(H^{\infty})$, which concludes the theorem.\\
\end{proof}
    
    Theorem \ref{Thereom characterization of Wg with Tg(Hoo)} implies that $T_g(H^{\infty})\subset \BMOA$ when $g\in\BMOA$. This was first shown in \cite{aleman1995integral}. 
    Due to Theorem \ref{Thereom characterization of Wg with Tg(Hoo)}, the spaces $W^p_g$ are independent of the exponent $p$ for $1\leq p<\infty$. Consequently, we will omit the superscript $p$ and simply write $W_g$.

\begin{proof}[Proof of Theorem \ref{T:equivalence meromorphic opttimal}]
 We first observe that $(T_{g_1},H^p)=(T_{g_2},H^p)$ if and only if $W_{g_1}=W_{g_2}$. 
Suppose $(T_{g_1},H^p)=(T_{g_2},H^p)$. Then for any  $f \in W_{g_1}$, we have $$(T_{g_1},H^p)\subseteq (T_f,H^p).$$
Since $(T_{g_2},H^p)=(T_{g_1},H^p)$, it follows that
$$(T_{g_2},H^p)\subseteq (T_f,H^p),$$ implying $f \in W_{g_2}$. Thus, $ W_{g_1}\subseteq W_{g_2}$. The reverse inclusion is followed by symmetry.

Suppose now that $ W_{g_1}= W_{g_2}$. Then $g_1 \in W_{g_1}=W_{g_2}$, so
$$(T_{g_2},H^p)\subseteq (T_{g_1},H^p)\,.$$ Similarly, $g_2 \in W_{g_1}$, so $(T_{g_1},H^p)\subseteq (T_{g_2},H^p)$. Hence $(T_{g_1},H^p)=(T_{g_2},H^p)$.
\par The desired conclusion in $b)$ follows, taking into consideration that for every $g\in\BMOA$ $W_{g}=T_g(H^{\infty})+\m{C}$, see Theorem \ref{Thereom characterization of Wg with Tg(Hoo)}. The functions $k_1,k_2\in H^{\infty}$ satisfy simultaneously the equations
$$\begin{cases}
    g_1'=k_1 g_2'\\
    g_2'=k_2g_1'
\end{cases}\, .$$
Therefore, it follows that $k_2=1/{k_1}\,.$\\
\end{proof}
 In the following proofs, we shall exploit the equivalent description of the $H^p$ norm of a function with a quantity that contains its derivative. Specifically, the Paley-Littlewood $G$-function is given by
$$
G(f)(t) = \left( \int_0^1 |f'(re^{it})|^2 (1 - r) \, dr \right)^{\frac{1}{2}} \quad \text{ for every } t \in \mathbb{R}\, .
$$
There exist two constants $C_1, C_2 > 0$, depending only on  $p$, such that
\begin{equation}\label{Radial function hp norm with derivatives}
C_1 \|f\|_{H^p}^p \leq |f(0)|^p + \int_{0}^{2\pi} G^p(f)(t) \, \frac{dt}{2\pi} \leq C_2 \|f\|_{H^p}^p\, ,
\end{equation}
for any $f$ analytic in $\m{D}$. See \cite[Chapter XIV, Theorem 3.19]{zygmund2002trigonometric} for the proof of \eqref{Radial function hp norm with derivatives}.
\begin{proof}[Proof of Theorem \ref{Theorem properties of Wg}]
(a)  Since $W_g=T_g(H^\infty)$ and the operator $T_g$ is bounded from $H^\infty$ to $\BMOA$, 
 the graph of $T_g$, which is precisely $W_g$, is a closed subspace of $H^\infty\times\BMOA$\CB.  Therefore, $W_g$ is a Banach space under $\|\cdot\|_{W_g}$, which is exactly the graph norm.
 
Moreover, for a fixed $z\in\m{D}$, since the point evaluations functionals are bounded both in $H^{\infty}$ and $\BMOA$, we have that
$$
|f(z)|\leq \lambda_z \|f\|_{W_g}\,  ,
$$
where $\lambda_z>0$ is a constant depending only on $z$. This shows that point evaluations are bounded functionals in $W_g$.

\par (b) The right inclusion is clear. For the left inclusion, we consider
\begin{equation}\label{equation Proposition intersection 1}
f \in \bigcap_{h \in W_g\colon  (T_h,H^p)\supsetneq (T_g,H^p)} (T_h,H^p)\, .
\end{equation}
First of all, we need to show that $fg'\in \Hol(\m{D})$. We assume the contrary, that is, we suppose that  $z_0\in\m{D}$ is a pole of $fg'$ of order $\lambda_1>0$. As $f$ satisfies \eqref{equation Proposition intersection 1}, we consider $h_1\in W_g$, such that $(T_{h_1},H^p)\supsetneq (T_g,H^p).$ Due to Theorem \ref{Thereom characterization of Wg with Tg(Hoo)}, there exists a function $q_{h_1}\in H^{\infty}$, such that $h'_1=g' q_{h_1}$. Since $fh'_1\in \Hol(\m{D})$, $q_{h_1}$ must vanish at $z_0$ of order $\lambda_2$ with $\lambda_2\geq \lambda_1$. Set
$$
q_{h_2}(z)=\frac{q_{h_1}(z)}{(z-z_0)^{\lambda_2}}(1-z)^{\frac{1}{p}},\qquad z\in\m{D}\, ,
$$
and let $h_2= T_g(q_{h_2})$. The function $q_{h_2}$ belongs to $H^\infty$, does not vanish at $z_0$ and, due to Theorem \ref{Thereom characterization of Wg with Tg(Hoo)},  $h_2\in W_g$. Moreover, $(T_{h_2},H^p)\supsetneq (T_g,H^p)$ since

$$ 
F(z)=\frac{1}{g'(z)}\frac{1}{(1-z)^{\frac{1}{p}+1}} \in (T_{h_2},H^p)\setminus (T_g,H^p)\, .
$$
Indeed
$$
T_g(F)(z)=\int_0^z \frac{1}{(1-\xi)^{\frac{1}{p}+1}}d\xi=p \left(\frac{1}{(1-z)^{\frac{1}{p}}}-1\right) \notin H^p\, ,
$$
and
$$
T_{h_2}(F)(z)=\int_{0}^zq_{h_2}(\xi)\frac{1}{(1-\xi)^{1/p+1}}d\xi=\int_{0}^z \frac{q_{h_1}(\xi)}{(\xi-z_0)^{\lambda_2}}\frac{1}{1-\xi}d\xi \in H^p\, ,
$$
since $q_{h_1}(z)/(z-z_0)^\lambda \in H^\infty$ and $-\log(1-z) \in H^p$.
\CB
As $f$ satisfies \eqref{equation Proposition intersection 1}, the function $fh_2'\in \Hol(\m{D})$, which is impossible, since 
$
fh'_2=fg'q_{h_2}
$
has a pole at $z_0$.

\par Now we prove that $T_g(f)\in H^p$. 
We chose $q_{1}(z)=(1+z)^{1/p}$ and $q_{2}(z)=(1-z)^{1/p}$.
Because of Theorem \ref{Thereom characterization of Wg with Tg(Hoo)}, $T_g(q_i) \in W_g$ and, as shown in the first part of the proposition, $(T_{T_g(q_i)},H^p)\supsetneq (T_g,H^p)$.
By applying \eqref{Radial function hp norm with derivatives}, we have that
\begin{align*}
    \|T_g(f)\|^p_{H^p}\lesssim &  \int_0^{2\pi}\left( \int_{0}^{1} (g'(re^{i\theta})f(re^{i\theta}))^2\left( 1-r \right)dr\right)^{p/2}d\theta\\
    =&  \int_{-\pi/2}^{\pi/2}\left( \int_{0}^{1} (g'(re^{i\theta})f(re^{i\theta}))^2\left( 1-r \right)dr\right)^{p/2}d\theta\\
    &\qquad \qquad+ \int_{\pi/2}^{3\pi/2}\left( \int_{0}^{1} (g'(re^{i\theta})f(re^{i\theta}))^2\left( 1-r \right)dr\right)^{p/2}d\theta\, .
\end{align*}
For the first integral, we have that
\begin{align*}
      &\int_{-\pi/2}^{\pi/2}\left( \int_{0}^{1} (g'(re^{i\theta})f(re^{i\theta}))^2\left( 1-r \right)dr\right)^{p/2}d\theta\\ 
      & \quad \leq \frac{1}{\inf_{z \in \mathbb D\cap \{\operatorname{Re}z >0\}}|q_1(z)|^p} \int_{-\pi/2}^{\pi/2}\left( \int_{0}^{1} (q_1(re^{i\theta}) g'(re^{i\theta})f(re^{i\theta}))^2\left( 1-r \right)dr\right)^{p/2}d\theta\\
    &\quad \lesssim  \|T_{T_g(q_1)}(f)\|_{H^p}^p<\infty\, ,
\end{align*}
and analogously for the second integral by considering the function $q_2$. This concludes the proof of the statement.

\par (c) 
    Let $\phi\in H^{\infty}$ and $T_\phi$ be a bounded operator on $W_g$. We prove that $\phi f\in W_g$ for every $f\in W_g$. Due to the Theorem \ref{Thereom characterization of Wg with Tg(Hoo)}, it is enough to prove that
    $
    (\phi f)'/g'\in H^{\infty}
    $.
Indeed
    \begin{align}\label{eqq3}
      \frac{(\phi(z) f(z))'}{g'(z)}&=\frac{\phi'(z) f(z)}{g'(z)}+\frac{\phi(z)f'(z)}{g'(z)} =\frac{T_\phi(f)'(z)}{g'(z)}+\phi(z)\frac{f'(z)}{g'(z)}\, .
    \end{align}
   Since $f'/g',T_\phi(f)'/g'\in  H^{\infty}$, we have that 
    $$
    \frac{(\phi f)'}{g'}\in H^{\infty}\, .
    $$
    \par 
     For the other implication, let $\phi$ be an analytic function which is a bounded multiplier of $W_g$.
    \CB
    Let $z\in\mathbb{D}$ and $\lambda_z$ the point evaluation functional in $W_g$. Since $W_g$ is a closed subspace of $\BMOA$ which contains the constant functions, $0\neq\|\lambda_z\|<\infty$. As already done in proof of Theorem \ref{Theorem properties of meromorphic optimal domain}.c, we obtain that
    $$
    |\phi(z)|\leq \|M_\phi\|\frac{\|\lambda_z\|}{\|\lambda_z\|}=\|M_\phi\|\, ,
    $$
    which implies that $\phi\in H^{\infty}$. Finally, we show that $T_\phi$ acts boundedly on $W_g$.
    Indeed, if $f \in W_g$, then
    $$
   \frac{\left(T_{\phi}(f)\right)'}{g'}=\frac{\phi'f}{g'}=\frac{(\phi f)'}{g'}-\frac{\phi f'}{g'}\, . 
    $$
    The first term is in $H^\infty$ since $\phi$ is a bounded multiplier of $W_g$. The second term belongs to $H^\infty$ since it is the product of two terms in $H^\infty$. So, $T_\phi(f)\in W_g$ for every $f\in W_g$ and consequently, a standard application of the closed graph theorem, implies the desired conclusion.\\
    \CB

\end{proof}

The first item shows that the space $W_g$ with $\|\cdot\|_{W_g}$ is a Banach space strictly contained in $\BMOA$.
A crucial question is whether $W_g$ is a closed subspace of $\BMOA$ under the $\|\cdot\|_{*}$ norm.
The authors have established a connection between this problem and the boundedness from below of the generalized Volterra operator $T_g$ acting from $H^\infty$ to $\BMOA$. Furthermore, they have demonstrated that, when $g$ is univalent, $\left( W_g,\|\cdot \|_{*}\right)$ is never closed, \cite{Preprint}.

We also highlight that the strength of item (b) lies in the fact that we choose functions with a strictly larger optimal domain than $(T_g,H^p)$: if we had considered the intersection of all $h\in W_g$, the result would have been trivial since $g\in W_g$ .

\vspace{22 pt}

\section{The space \texorpdfstring{$V_g^p$}{} and equivalence of holomorphic optimal domains}

In this section, we give a brief insight into the space $V_g^p$ and investigate the conditions under which the holomorphic optimal domains of two $\BMOA$ functions, $g_1$ and $g_2$, coincide.

A complete characterization of $V^p_g$ for arbitrary $g \in \BMOA$ remains elusive, and establishing an analog of Theorem \ref{Thereom characterization of Wg with Tg(Hoo)} for $V^p_g$ \CB poses significant challenges.
First of all, we note that $V^p_g$ is independent of the choice of $p$ when $g$ is locally univalent.
\CB
\begin{prop}\label{Proposition characterisation Vg g loc univ}
 Let $g\in\BMOA\cap U_{\loc}(\m{D}) $ and $1\leq p<\infty$. $h\in V_g^p$ if and only if there exists a function $f\in T_g(H^{\infty})$ such that $h=f+h(0)$. In other words,
 $$V_g^p=T_g(H^{\infty})+\m{C}\,.$$   
\end{prop}
\begin{proof}
The proof follows the lines of Theorem \ref{Thereom characterization of Wg with Tg(Hoo)}, once we note that, when $g \in U_{\loc}(\m{D})$,
\begin{equation}\label{ppp}
    [T_g,H^p]=\left\lbrace f \in \text{Hol}(\m{D})\colon f(z)=\frac{h'(z)}{g'(z)} \text{ for some }h \in H^p\right\rbrace\,  .
\end{equation}
\end{proof}

When $g\in U_{\text{loc}}(\m{D})$, Proposition \ref{Proposition characterisation Vg g loc univ}, implies that $V_g^p=T_g(H^{\infty})+\m{C}$, providing that $V_g^p$ is independent of $p$, however, we do not have this characterization when $g$ is not a locally univalent function. 
Indeed, let $g(z)=z^2$, which is not locally univalent. \CB It is easily checked that $[T_{z^2},H^p]=[T_z,H^p]$, although this result will be immediate taking into account Proposition \ref{Theorem optimal Domain polynomial}.  Then, due to Proposition \ref{Proposition characterisation Vg g loc univ}, $V_{z^2}^p=T_z(H^{\infty})+\m{C}$. However, 
$$
z+z^2=T_z(1+t/2)\in V^p_{z^2}\text{ but } z+z^2\notin T_{z^2}(H^{\infty})\, 
$$
since the derivative of $z+z^2$ does not vanish at $0$.

Our most significant achievement in characterizing $V^p_g$ for general $g \in \BMOA$ is presented in the following proposition.
\CB

\begin{prop}\label{proposition V_g when g in nevanlinna}
Let $1\leq p<\infty$ and $g\in \BMOA$. Then 
$$
T_g(H^\infty)+\mathbb{C}\subseteq V_g^p
$$
Moreover, \CB if there exists a function $u\in \BMOA$ and $B \in H^\infty$ such that $g'=Bu'$, then
 $$ V_g^p \subseteq  V^p_u\,.$$
\end{prop}
\begin{proof}
    Let $h=T_{g}(q)+h(0)$ for some $q\in H^\infty$. For the first inclusion, it suffices to prove that 
$T_h(f)\in H^p$ for every $f\in[T_g, H^p]$.
By \cite[Theorem 4]{Bellavitaoptimaldom}, we have that $fq\in [T_g,H^p]$ and, consequently,
\begin{align*}
T_{h}(f)=T_{g}(qf)\in H^p\, .
\end{align*}
For the second inclusion, we use the representation $g'=Bu'$ and we show that $V_g^p\subseteq V_u^p$. To this end, we observe that for every $f \in [T_u,H^p]$,
$$
T_g(f)=T_u(Bf) \in H^p\, ,
$$
due to \cite[Theorem 4]{Bellavitaoptimaldom}. It follows that $[T_u,H^p]\subseteq [T_g,H^p]$, which proves the statement.\\
\end{proof}

In the rest of the article, \CB we answer for which functions $g_i$ the corresponding holomorphic optimal domains coincide assuming extra conditions on $g_i \in \BMOA$. We first present an equivalent description of the problem, when $g_2'/g_1'$ is analytic on the unit disc.
\begin{prop}\label{general optimaldomain question}
    Let $1\leq p<\infty$ and $g_1,g_2$ in $\BMOA$ such that $g'_2/g'_1=h \in \Hol(\mathbb D)$. $[T_{g_1},H^p]=[T_{g_2},H^p]$ if and only if $h \in H^\infty(\mathbb D)$ and $M_h$ is a multiplier of  $[T_{g_1},H^p]$ bounded from below. 
\end{prop}
\begin{proof}
    First of all, let us suppose that $[T_{g_1},H^p]=[T_{g_2},H^p]$. Then, because of the open mapping theorem, for every $f \in [T_{g_1},H^p]$
    $$
    \|f\|_{[T_{g_1},H^p]}\sim \|f\|_{[T_{g_2},H^p]}\, .
    $$
    Moreover,
$$
\|hf\|_{[T_{g_1},H^p]}=\|T_{g_1}(hf)\|_{H^p}=\|T_{g_2}(f)\|_{H^p}=\|f\|_{[T_{g_2},H^p]}\lesssim \|f\|_{[T_{g_1},H^p]}\, .
$$
Consequently $h \in \mathcal{M}([T_{g_1},H^p])$ and, because of \cite[Theorem 4]{Bellavitaoptimaldom}, $h \in H^\infty(\mathbb D)$.
On the other hand, for every $f \in [T_{g_1},H^p]$,
\begin{align*}
    \|hf\|_{[T_{g_1},H^p]}&=\|T_{g_2}(f)\|_{H^p}=\|f\|_{[T_{g_2},H^p]}\gtrsim  \|f\|_{[T_{g_1},H^p]}\, .
\end{align*}
 The sufficiency can be proved by using the same inequalities.\\
\end{proof}
Of course, this result does not provide any further insight on the problem since the multipliers of  $[T_{g},H^p]$ bounded from below have not been characterized. \CB For this reason, we study the equivalence of two holomorphic optimal domains when the functions $g_i$ satisfy some extra natural conditions.
\subsection{Locally univalent}
When the function $g$ is locally univalent,
we have an explicit description for the space $V_g^p$.

As shown in Theorem \ref{T:equivalence meromorphic opttimal}, in this situation we establish when two holomorphic optimal domains contain the same functions.
\begin{corollary}\label{prop univalent}
 For $1\leq p<\infty$, let $g_1,g_2\in\BMOA$, non-constant. The following are equivalent;
 \begin{itemize}
     \item[a)] $[T_{g_1},H^p]=[T_{g_2},H^p]$;
     \item[b)] There exist two functions $k_1,k_2\in H^{\infty}$ such that 
     $$
g_1=T_{g_2}(k_1)-g_1(0) \text{ and } g_2=T_{g_1}(k_2)-g_2(0)\, .
$$
 \end{itemize}
 Moreover we have that $k_2=\frac{1}{k_1}$\,.  
\end{corollary}

\subsection{Polynomials}
 The second case, we consider is where the symbols $g_i$ are analytic polynomials. 
 We need the following auxiliary lemma.

\begin{Lemm}\label{Division by polynomial with zeors in D}
Let $1\leq p<\infty $ and $g \in \BMOA$. If $g'=Qu'$ with $u \in \BMOA$ and $Q$ an analytic polynomial with no zeros in $\m{T}$, then $[T_g,H^p]=[T_u,H^p]$.
\end{Lemm}
\begin{proof}
Let $ f \in [T_u,H^p]$. Because of \cite[Theorem 4]{Bellavitaoptimaldom}, $Qf \in [T_u,H^p]$. Therefore 
$$
T_g(f)=T_u(Qf) \in H^p
$$
and $[T_u,H^p]\subseteq[T_g,H^p] $.

To prove the reverse implication, we consider $f$ belonging to $[T_g,H^p]$. If the polynomial $Q$ does not have any root in $\overline{\m{D}}$, then $Q$ is a function bounded from below on $\m{D}$, and the inclusion can be easily checked, taking into account \eqref{Radial function hp norm with derivatives}. To this end, we shall work with $Q$ having roots inside the unit disc and the proof shall be done by induction on the degree of $Q$. When $n=1$, let $Q(z)=z-z_0$ for some $z_0\in\m{D}$. Set 
$$h(z)=\int_{z_0}^zf(\zeta)g'(\zeta)\,d\zeta\qquad z\in\m{D}\,.$$
Then, $h \in H^p$, with $h(z_0)=h'(z_0)=0$ and $g'f=h'$. Consequently
\begin{align*}
    T_u(f)(z)&=\int_0^z u'(\xi)f(\xi)d\xi=\int_0^z u'(\xi)\frac{h'(\xi)}{g'(\xi)}d\xi=\int_0^z \frac{h'(\xi)}{(\xi-z_0)}d\xi\\
    &=\frac{h(z)}{z-z_0}+\frac{h(0)}{z_0}+\int_{0}^z\frac{h(\xi)}{(\xi-z_0)^2}d\xi\, .
\end{align*}
As $z_0\in\m{D}$, the function $h(z)/(z-z_0)\in H^p$. Moreover, because of \cite[Theorem 3.11]{duren1970theory}, the function $\int_0^z {h(\xi)}/{(\xi-z_0)^2}d\xi$ belongs to the disk algebra and consequently to $H^p$. Due to this reasoning, $f \in [T_u,H^p]$. 
\par If the degree of $Q$ is $n$, let $\mathcal{Z}(Q)=\{z_0,\dots,z_{n-1}\}$ be the zero set of $Q$, each point taken according to its multiplicity. We set 
$$g_{n-1}(z)=\int_{0}^{z}(\zeta-z_0)\dots(\zeta-z_{n-2})u'(\zeta)\,d\zeta\qquad z\in\m{D}\,.$$
Then, $g_{n-1}'(z)=(z-z_0)\cdots(z-z_{n-2})u'(z)$ for $z\in\m{D}$ and by the induction hypothesis, we have that
\begin{equation}\label{equation1lemma}
    [T_{g_{n-1}},H^p]=[T_u,H^p]\,.
\end{equation}
However, the formula of $g'$ implies that $g'=(z-z_{n-1})g_{n-1}'$ and consequently the first part of the proof implies that
\begin{equation}\label{equation2lemma}
    [T_g,H^p]=[T_{g_{n-1}},H^p]\,.
\end{equation}
By combining \eqref{equation1lemma} and \eqref{equation2lemma}, we have the desired conclusion.\\
\end{proof}
  \begin{prop}\label{Theorem optimal Domain polynomial}
    Let $1\leq p< \infty$ and $g$ be an analytic polynomial. Let $\lambda_j\in\overline{\m{D}}$ be the points where $g'(\lambda_j)=0.$ The following hold:
    \begin{itemize}
        \item[(a)] If $|\lambda_j|\neq 1$ for every $j$, $[T_g,H^p]=[T_z,H^p]$.
        \item[(b)] If there exists an index $j$ such that $|\lambda_j|=1$, then $[T_z,H^p]\subsetneq [T_g,H^p]$.
    \end{itemize}
\end{prop}
\begin{proof}
\par The part (a) follows from Lemma \ref{Division by polynomial with zeors in D}. 
For the part (b), we prove the proposition assuming that $\lambda_1=1$. The general case is followed by a rotation argument. 
Let $g'(z)=(1-z)p(z)$, where $p$ is an analytic polynomial. 
Since $[T_z,H^p]=[T_{z^n},H^p]$ because of Lemma \ref{Division by polynomial with zeors in D}, it is also clear that 
$$
[T_z,H^p]\subset [T_g,H^p]\, .
$$
On the other hand, we notice that the function
$(1-z)^{-1-1/p} \in [T_{(1-z)^2},H^p]$, but it does not belong to $[T_z,H^p]$. Indeed, 
$$
T_z\left(\frac{1}{(1-z)^{\frac{1}{p}+1}}\right)(\lambda)= \frac{1}{p} \left(\frac{1}{(1-\lambda)^{{\frac{1}{p}}}}-1\right)\not\in H^p\qquad \text{ for }\lambda\in\m{D}\, .
$$
On the other hand, for $\lambda\in\m{D}$,
$$
T_{(1-z)^2}\left(\frac{1}{(1-z)^{\frac{1}{p}+1}}\right)(\lambda)= \begin{cases}
    -\frac{p}{p-1} \left(1-\lambda)^{1-\frac{1}{p}}-1\right),& p\neq1\\
    -2\log\frac{1}{1-\lambda},& p=1
\end{cases}\, 
$$
which belongs to $H^p$ for any $1\leq p <\infty$. Therefore, $[T_z,H^p] \subsetneq [T_{(1-z)^2},H^p]$. Finally, $[T_{(1-z)^2},H^p]\subseteq [T_{g},H^p]$. For, if $f\in [T_{(1-z)^2},H^p]$, then, since $p$ is a polynomial and consequently a bounded function, we have that $fp\in [T_{(1-z)^2},H^p]$ as a result of \cite[Theorem 4]{Bellavitaoptimaldom}. Moreover,
$$T_{g}(f)(z)=\frac{1}{2}T_{(1-z)^2}(fp)(z)\qquad z\in\m{D},$$
which proves that $f\in [T_{g},H^p]$, as well. Combining the last two inclusions, the proof is complete.\\
\end{proof}

\subsection{The Hilbert case}
 In this situation, we restrict ourselves to the holomorphic optimal domains when $p=2$: we look for cases when given $g_1,g_2 \in \BMOA$, $[T_{g_1},H^2]=[T_{g_2},H^2]$. We recall from \cite{Bellavitaoptimaldom} that 
$$
[T_g,H^2]=\biggr\{f\in \Hol(\m{D})\colon \int_{\m{D}}|f(z)|^2|g'(z)|^2(1-|z|^2)\,dm(z)<\infty\biggl\}=:A^2(|g'(z)|^2(1-|z|^2))
$$
where $dm(z)$ is the normalized Lebesgue area measure of the unit disc. This case is worth studying since Proposition \ref{general optimaldomain question} implies the following result.
\begin{corollary}\label{p=2 optimaldomain question}
    Let $g_1,g_2$ in $\BMOA$ such that $g'_2/g'_1=h \in \Hol(\mathbb D)$. $[T_{g_1},H^2]=[T_{g_2},H^2]$ if and only if $h \in H^\infty(\mathbb D)$ and $M_h$ is a multiplier of the Bergman space $A^2\left( |g_1'(\xi)|^2\left( 1-|\xi|^2\right)\right)$ bounded from below. 
\end{corollary}
The space $A^2(|g'(\xi)|^2(1-|z|^2))$ is a non-radial weighted Bergman space of analytic functions. Due to the limited research on such spaces, as noted in \cite{pelaez2014weightedmemoir} and \cite{pelaez2015smallweighted}, characterizing the multipliers bounded from below and, consequently, determining when the holomorphic optimal domains $[T_{g_1},H^2], [T_{g_2},H^2]$ are equal (as per Corollary \ref{p=2 optimaldomain question}) is challenging.

To address this problem, we focus on a specific class of functions $g$'s: we consider functions $g$'s whose derivatives can be written as the product of a Blaschke product $B$ and a non-vanishing function $u'$.

\begin{prop}\label{corollary of horowitz}
    Let $g\in\BMOA$ such that $g'=Bu'$, where $B$ is a finite product of interpolating Blaschke products, $u \in \BMOA\cap U_{\text{loc}}(\m{D})$. Then, $[T_g,H^2]=[T_u,H^2]$.
\end{prop}
\begin{proof}
We first show that
\begin{equation}\label{equation 1prop 5.8}
    [T_g,H^2]=\{f\in \Hol(\m{D})\colon f(z)=\frac{h(z)}{u'(z)}\, \, \text{ where }  \quad Bh\in A^2(1-|z|^2)\}\, .
\end{equation}
Let $f\in [T_g,H^2]$. Then, $T_g(f)\in H^2$. Set $H=T_g(f)'$. Then, $H'\in A^2(1-|z|^2)$ and moreover
$$f(z)=\frac{H(z)}{B(z)}\frac{1}{u'(z)}\qquad z\in\m{D}\, .
$$
Setting $h=H/B$ provides the desired inclusion. On the other hand, let $f$ belong to the right-hand side of \eqref{equation 1prop 5.8}. Taking into account the factorization of $g'$, standard calculation shows that $T_g(f)\in H^2$, hence $f\in [T_g,H^2]$.
\par Since the Blaschke product $B$ is a finite product of interpolating Blaschke products, the result of C. Horowitz \cite{Horowitz1977} implies that $B$ is a universal divisor of $A^2(1-|z|^2)$. Consequently
\begin{equation*}\label{equation horowitz result}
    [T_g,H^2]=\{f\in\Hol(\m{D})\colon f(z)=\frac{h(z)}{u'(z)}\quad h\in A^2(1-|z|^2)\}=[T_u,H^2]\, ,
\end{equation*}
where the second equality follows from \eqref{ppp}. The proof is complete.\\
\end{proof}

\vspace{22 pt}
\section{Acknowledgements}
The authors are grateful to Professor A. Siskakis for stimulating discussions on the topic of this article and for his insightful comments on the first draft of this article.

\par
The first author is a member of Gruppo Nazionale per l’Analisi Matematica, la Probabilit\`a e le loro Applicazioni (GNAMPA) of Istituto Nazionale di Alta Matematica (INdAM) and he was partially supported by PID2021-123405NB-I00 by the Ministerio de Ciencia e Innovaci\'on and by the Departament de Recerca i Universitats, grant 2021 SGR 00087.

The third and the fourth authors were partially supported by the Hellenic Foundation for Research and Innovation (H.F.R.I.) under the `2nd Call for H.F.R.I. Research Projects to support Faculty Members \& Researchers' (Project Number: 4662).

\vspace{22 pt}
\bibliography{bibliography}
\bibliographystyle{plain}
\vspace{22 pt}

\end{document}